\documentclass[12pt,letterpaper]{amsart}   

\usepackage{amsmath,amsthm,amscd,amssymb,array,transfig}

\newtheorem{theorem}{Theorem}
\newtheorem{lemma}{Lemma}
\theoremstyle{remark}
\newtheorem*{remark*}{Remark}
\newcommand{\proj}{\mathbf{P}}
\newcommand{\Q}{\mathbb{Q}}
\newcommand{\Z}{\mathbb{Z}}

\renewcommand{\P}{\mathbf{P}}

\renewcommand{\emptyset}{\varnothing}
\newcommand{\Hbar}{\overline{\mathcal{H}}}
\newcommand{\Mbar}{\overline{\mathcal{M}}}

\newcommand{\etc}{\,\ldots}

\newcommand{\into}{\hookrightarrow}

\newcolumntype{C}{>{$}c<{$}}    
\newcolumntype{L}{>{$}l<{$}}    
\newcolumntype{R}{>{$}r<{$}}    
\newcommand{\barr}{\overline}

\newcommand{\sy}{\mathbb{S}}
\newcommand{\G}{\mathbf{G}}

\newcommand{\mgn}{\overline{\mathcal{M}}_{g,n}}
\newcommand{\mgnX}{\overline{\mathcal{M}}_{g,n}(X,\beta)}
\newcommand{\mtt}{\overline{\mathcal{M}}_{2,3}}
\newcommand{\scup}{\mathbin{\text{\scriptsize$\cup$}}}
\newcommand{\scap}{\mathbin{\text{\scriptsize$\cap$}}}
\newcommand{\lan}{\langle}
\newcommand{\ran}{\rangle}
\newcommand{\virt}{{\textup{virt}}}
\newcommand{\ns}{{\textup{ns}}}
\newcommand{\sing}{{\textup{sing}}}
\newcommand{\red}{{\textup{red}}}
\newcommand{\ev}{{\textup{ev}}}
\newcommand{\sym}{{\textup{sym}}}
\newcommand{\flag}{{\textup{flag}}}
\newcommand{\disj}{{\scriptscriptstyle\sqcup}}
\newcommand{\R}{N^{\scriptscriptstyle(0)}}
\newcommand{\E}{N^{\scriptscriptstyle(1)}}
\newcommand{\N}{N^{\scriptscriptstyle(2)}}
\newcommand{\HH}{H^{\scriptscriptstyle(2)}}
\newcommand{\PP}{P^{\scriptscriptstyle(2)}}

\newcommand{\open}{\bigg[\raisebox{-2ex}}
\newcommand{\close}{\hskip-3pt \bigg]}

\begin{document}
\title{A Descendent Relation In Genus 2}
\author{P. Belorousski and R. Pandharipande}
\address{Department of Mathematics, University of Chicago, Chicago, IL
  60637, USA}
\email{pavel@math.uchicago.edu \\rahul@math.uchicago.edu}
\maketitle

\pagestyle{plain}
\setcounter{section}{-1}
\section{\bf Introduction}

Let $\mgn$ be the moduli space of Deligne-Mumford stable, $n$-pointed, 
genus $g$, complex algebraic curves. There is an algebraic
stratification of $\mgn$ by the underlying topological type of the
pointed curve. These strata are naturally indexed by stable,
$n$-pointed, genus $g$ dual graphs. A relation among the cycle classes 
(or homological classes) of the closures of these strata directly
yields  differential equations satisfied by generating functions of
Gromov-Witten invariants  of algebraic varieties. The translation from 
a relation to differential equations is obtained by the
splitting axiom of Gromov-Witten theory (\cite{ruantian1},
\cite{kontsmanin1}, \cite{behrmanin}). In genus 0, 
all strata relations are obtained from the basic linear equivalence of
the three boundary strata in $\overline{\mathcal{M}}_{0,4}$
(\cite{kontsmanin3}). The corresponding differential equation is the 
Witten-Dijkgraaf-Verlinde-Verlinde equation. In genus 1, Getzler
has found a codimension 2 relation in $\overline{\mathcal{M}}_{1,4}$
(\cite{getzler1}).   
The resulting differential equation has been used to calculate
elliptic Gromov-Witten invariants (\cite{getzler1})
and to prove a genus 1 prediction of the Virasoro conjecture for
$\proj^2$ (\cite{rahul}). Getzler's equation has been studied in the
context of semi-simple Frobenius manifolds in \cite{kabkim} and
\cite{dubzha}.

Let $X$ be a nonsingular complex projective variety. Let $\mgnX$
be the moduli space of stable maps representing the class 
$\beta \in H_2(X, \Z)$. 
The {\em descendent invariants} are the integrals
\begin{equation}
\label{grav}
\lan (\gamma_1\psi_1^{\alpha_1})\cdots
(\gamma_n\psi_n^{\alpha_n})\ran_{g,\beta} = 
\underset{[\mgnX]^\virt}{\int} \ev^*_1(\gamma_1)\scup
\psi_1^{\alpha_1} \scup 
\cdots \scup \ev_n^*(\gamma_n)\scup \psi_n^{\alpha_n},
\end{equation}
where $\gamma_i\in H^*(X, \Q)$ and $\psi_i$ is the first Chern class
of the cotangent line on (the stack) $\mgnX$ corresponding to
the $i$-th marked point. These invariants play a central
role in Gromov-Witten theory. For example,  they arise 
naturally in expressions for flat sections in the Dubrovin
formalism (\cite{dub}, \cite{giv}), in the virtual normal bundle terms
in torus localization formulas (\cite{konts}, \cite{giv},
\cite{grabpand1}), and in the Virasoro conjecture \cite{ehx}. 
A geometric interpretation of certain low genus descendent invariants
of $\proj^2$ in terms of the classical characteristic numbers of plane 
curves is given in \cite{grabpand2}. 
Some foundational issues concerning descendent integrals
are treated in \cite{ruantian2}, \cite{kontsmanin2}, and
\cite{getzler2}.
 
In genus 0 and 1, the classes of the cotangent lines on $\mgn$ may be 
expressed as sums of boundary divisor classes. Such expressions 
yield {topological recursion relations} among the descendent
invariants and may be used to prove that in genus 0 and 1 the
descendents are determined by Gromov-Witten invariants. In higher
genera, divisorial topological recursion relations do not exist. In
genus 2, Getzler has determined weaker topological recursion relations
from boundary expressions for $\psi_1^2$ and $\psi_1\psi_2$
(\cite{getzler2}).  

The  strata in $\mgn$ are themselves products of moduli spaces of  
pointed curves (quotiented by symmetry groups).
Therefore, cotangent line classes are defined
on the strata.  It should be noted that 
not all of these classes  are restrictions of cotangent
lines on $\mgn$.
A {\em descendent stratum class} in $A^*(\mgn)$ is 
the push-forward from a stratum $S$ to $\mgn$ of
a monomial in the cotangent lines of $S$.
Relations among the descendent stratum classes yield differential
equations for the generating functions  of descendent invariants
of algebraic varieties. In physics, this generating function $F_X$ for
a nonsingular projective variety $X$ is called the full gravitational
potential function. The Virasoro conjecture states that 
${\rm exp}(F_X)$ is annihilated by an explicit  representation
$\rho_X$ of the affine Virasoro algebra in an algebra of differential
operators (\cite{ehx}).  

In this paper we present a new genus $2$ relation among codimension 2
descendent stratum classes in $\mtt$. Getzler has computed
$h^4(\mtt)=44$ via a subtle method using mixed Hodge theory and
modular operads (\cite{getzler2}). The number of descendent stratum
classes in $A^2(\mtt)$ is $47$. Exactly 2 relations are obtained from
the basic genus 0 linear equivalence.
Therefore, there must exist a new relation, at least in homology.
We find an algebraic relation in $A^2(\mtt)$ via the admissible cover
technique introduced in \cite{rahul}.  
In fact, the relation lies in the $\mathbb{S}_3$-invariant subspace.
The resulting differential equations together with the
known topological recursion relations  are strong enough
to determine  the genus 2 descendent integrals for $\proj^2$.

The plan of the paper is as follows. 
The admissible double cover construction is reviewed in Section 1.
In Section 2, this construction is used to calculate the new relation
in $A^2(\mtt)$. Formulas expressing the cycle classes of Weierstrass
loci in the moduli spaces of pointed curves of genus 1 and 2 in terms
of descendent stratum classes are required for the computation of the
new relation. These formulas are also  obtained in Section 2 via the
space of admissible double covers. The application to the descendent
integrals of $\proj^2$ appears in Section 3.

Our greatest mathematical debt is to E. Getzler for informing us of
his homological computation.  His work provided the motivation for our
calculation. Thanks are also due to C. Faber and T. Graber for
related conversations. The second author was partially supported
by a National Science Foundation post-doctoral fellowship.

\section{\bf Admissible double covers}
\label{admi}
A genus $g$  {\em admissible cover} with $n$ marked points 
and $b$ branch points consists
of a morphism $\pi: C\to D$ of pointed curves
$$(C,\ P_1,\etc,P_n),\ \ (D,\ p_1,\etc,p_n,q_1,\etc,q_b)$$ 
satisfying the following conditions.
\begin{enumerate}
\item[(1)] $C$ is a connected, reduced, nodal curve of arithmetic genus $g$. 
\item[(2)] The markings $P_i$ lie in the nonsingular locus $C_\ns$.
\item[(3)] $\pi(P_i)=p_i$.
\item[(4)] $(D,\ p_1,\etc,p_n,q_1,\etc,q_b)$ is an $(n+b)$-pointed
  stable curve of genus $0$.
\item[(5)] $\pi^{-1}(D_{\ns}) = C_{\ns}$.
\item[(6)] $\pi^{-1}(D_{\sing})=C_{\sing}$.
\item[(7)] $\pi|_{E_{\ns}}$ is \'etale except over the points $q_j$ where
  $\pi$ is simply ramified. 
\item[(8)] If $x\in C_{\sing}$, then
\begin{enumerate} 
\item[(a)] $x \in C_1 \cap C_2$, where $C_1$ and $C_2$ are distinct
  components of $C$,
\item[(b)] $\pi(C_1)$ and $\pi(C_2)$ are distinct components of $D$,
\item[(c)] the ramification numbers at $x$ of the two morphisms
  $$\pi: C_1 \to \pi(C_1)\ \ \text{and}\ \ \pi:C_2 \to \pi(C_2)$$ are equal.
\end{enumerate}
\end{enumerate} 
These conditions imply that the map $\pi: C \to D$ is of uniform
degree $d$, where $$-2d+b=2g-2.$$
Let $\Hbar_{d,g,n}$ be the space of $n$-pointed genus $g$
admissible covers of $\proj^1$ branched at $b$ points.

Only the admissible double cover case in genus 1 and 2 will be
considered in this paper. The space
$\Hbar_{2,g,n}$ is an irreducible variety.
There are natural morphisms
\begin{equation}
\label{mmm}
\lambda:\Hbar_{2,g,n} \to \Mbar_{g,n}, \quad
\pi:\Hbar_{2,g,n} \to \Mbar_{0,n+2g+2}
\end{equation}
obtained from the domain and range of the admissible cover
respectively. In genus 1 and 2, $\lambda$ is clearly surjective.
The projection $\pi$ is a finite map to $\Mbar_{0,n+2g+2}$.
For each marking $i \in\{1,\etc,n\}$, there is a
natural $\Z/2\Z$-action given by switching
the sheet of the $i$-th marking of $E$.
These actions induce a product action in which the
diagonal $\Delta$ acts trivially. Define the group $\G$ by:
 $$\G= (\Z/2\Z)^n/\Delta.$$
The action of $\G$ on $\Hbar_{2,g,n}$ is generically free and commutes
with the projection $\pi$. Therefore, the quotient $\Hbar/\G$
naturally maps to $\Mbar_{0,n+2g+2}$. In fact, since the morphism
$$
\Hbar/\G \to \Mbar_{0,n+2g+2}
$$ 
is finite and birational, it is an isomorphism.

Spaces of admissible covers were defined in \cite{harmum}. The
methods there may be used to construct spaces of pointed admissible
covers. An alternative construction of the space of pointed admissible
covers via Kontsevich's space of stable maps is given in \cite{rahul}. 
A foundational treatment of the moduli problem of admissible
covers is developed in \cite{abrvist}.

Our method of obtaining a codimension 2 relation in $\Mbar_{2,3}$
is the following. Consider the diagram:
\begin{equation}
\label{pushpull}
\begin{CD}
\Hbar_{2,2,3} @>{\lambda}>> \Mbar_{2,3} \\ 
@V{\pi}VV \\
\Mbar_{0,3+6} \\ 
\end{CD}
\end{equation}
There are relations in $A^2(\overline{\mathcal{M}}_{0,3+6})$ among the
classes of codimension 2 strata. Such relations yield cycle relations in
$A^2(\mtt)$ via pull-back by $\pi$ and push-forward by $\lambda$.
The resulting cycles in $\mtt$
include Weierstrass loci. By further expressing these
Weierstrass loci in terms of descendent strata, a new nontrivial relation
is obtained.

\section{\bf The relation computation}
\label{relcomp}

\subsection{The new relation}
The main result of this section is the following theorem.
\begin{theorem} The codimension 2 descendent stratum classes in
  $\Mbar_{2,3}$ satisfy a nontrivial rational equivalence:
\begin{equation}
\label{maine}
\begin{split}\textstyle
&-2\open{\input{tr1.pictex}}\close + 2\open{\input{p.pictex}}\close + 3\open{\input{p1.pictex}}\close - 3\open{\input{p2.pictex}}\close \\[3pt]
&+\frac{2}{5}\open{\input{er3.pictex}}\close - \frac{6}{5}\open{\input{er2.pictex}}\close + \frac{12}{5}\open{\input{er1.pictex}}\close 
  -\frac{18}{5}\open{\input{e1r1.pictex}}\close - \frac{6}{5}\open{\input{ee.pictex}}\close \\[3pt]
&+\frac{9}{5}\open{\input{e1e.pictex}}\close - \frac{6}{5}\open{\input{ee1.pictex}}\close + \frac{1}{60}\open{\input{rn3.pictex}}\close
  -\frac{3}{20}\open{\input{rn2.pictex}}\close + \frac{3}{20}\open{\input{rn1.pictex}}\close \\[3pt]
&-\frac{1}{60}\open{\input{rn.pictex}}\close + \frac{1}{5}\open{\input{b0.pictex}}\close - \frac{3}{5}\open{\input{b1.pictex}}\close
  +\frac{1}{5}\open{\input{b2.pictex}}\close - \frac{1}{10}\open{\input{en.pictex}}\close - \frac{1}{10}\open{\input{en1.pictex}}\close = 0.
\end{split}
\end{equation}
\end{theorem}

We begin by explaining our notation.
All classes considered are {\em stack} fundamental
classes, and all Chow groups are taken with rational coefficients. 
A stratum $S$ is denoted by the topological type of the stable curve
corresponding to the generic moduli point of $S$. In the diagrams,
the geometric genera of the components are underlined.

A descendent stratum class in $\Mbar_{2,3}$ with {\em unassigned}
markings denotes the sum of descendent stratum classes over the 
$3!$ marking assignment choices. 
For example,
$$\open{\input{tr1.pictex}}\close = 2\open{\input{tr1ordc.pictex}}\close + 2\open{\input{tr1ordb.pictex}}\close + 2\open{\input{tr1ord.pictex}}\close.$$
The main advantage of this notation is the following:
an unsymmetrized equation among descendent strata
may be symmetrized by simply erasing the markings.
While our main relation (\ref{maine}) is symmetric,
most of the auxiliary formulas expressing classes of geometric loci 
in terms of descendent stratum classes are not. Since many of these
formulas are of independent interest, we present them in the 
unsymmetrized form. 

There are 3 strata in equation (\ref{maine}) with cotangent line
classes:
$$\open{\input{p.pictex}}\close, \open{\input{p1.pictex}}\close, \open{\input{p2.pictex}}\close.$$
The cotangent line class is always on the genus 2 component.
On the first and third
strata above, the cotangent line is taken at the node; on the
second stratum, it is taken at the marked point. 

The new relation will be found via an admissible double cover
construction. Let $\Mbar_{0,3+6}$ be the moduli space of stable genus
zero curves with the marking set $$\{p_1, p_2, p_3, b_1,\etc,b_6\}.$$  
The 9-pointed genus zero curve will  be the base of the admissible
cover: the points $p_i$ correspond to the images downstairs of the
marked points of the cover and the points $b_i$ correspond to the
branch points.  

\begin{remark*}
Let $\lambda$ and $\pi$ be the morphisms from diagram
(\ref{pushpull}). Let $\sy_6$ act on $\Mbar_{0,3+6}$ by permuting
the branch points $b_i$. Then the homomorphism 
$$\lambda_*\pi^*:A^2(\Mbar_{0,3+6})\to A^2(\Mbar_{2,3})$$
is $\sy_6$-invariant. 
\end{remark*}

\subsection{The first equation}
\label{firsteq}
Let $D$ denote the boundary divisor $D(p_1p_2|p_3b_1\ldots b_6)$ in
$\Mbar_{0,3+6}$. As a variety it is isomorphic to
$\Mbar_{0,A}$ with the marking set
$$A=\{\star,p_3,b_1,\etc,b_6\}.$$ Here $\star$ denotes the
node on the stable curve corresponding to the generic point of
$D$. Consider the following 4-point linear equivalence of divisors on 
$\Mbar_{0,A}$ obtained via boundary equivalence
on $\Mbar_{0,\{\star,p_3,b_1, b_2\}}$: 
$$ (\star p_3|b_1b_2)\sim (p_3b_1|\star b_2). $$
Let $\Lambda_D$ denote the corresponding relation in
$A^2(\Mbar_{0,3+6})$. Symmetrization with respect to the
natural $\sy_3$-action on $\Mbar_{0,3+6}$ permuting the points $\{p_1,
p_2, p_3\}$ yields the relation 
\begin{equation}
\label{wdvv1} 
\sum_{\sigma\in\sy_3} \sigma^* \Lambda_D
\end{equation}
in $A^2(\Mbar_{0,3+6})$.

\begin{lemma} 
\label{lemma1}
The application of $\lambda_*\pi^*$ to relation (\ref{wdvv1}) yields
(\/$6!$ times) the following rational equivalence in
$A^2(\Mbar_{2,3})$: 
\begin{equation}
\label{rel1}
\begin{split}
&\open{\input{tr1.pictex}}\close + 3\open{\input{capee.pictex}}\close + \frac{4}{3}\open{\input{d.pictex}}\close + \frac{2}{5}\open{\input{b4.pictex}}\close 
        + \frac{2}{5}\open{\input{ee1.pictex}}\close + \frac{2}{5}\open{\input{b0.pictex}}\close  \\[3pt] 
&+ \frac{1}{30}\open{\input{en1.pictex}}\close + \frac{1}{30}\open{\input{c5.pictex}}\close = \frac{1}{3}\open{\input{g.pictex}}\close 
        + \frac{2}{3}\open{\input{d2.pictex}}\close + \frac{1}{3}\open{\input{d1.pictex}}\close \\[3pt] 
&+ \frac{3}{5}\open{\input{b3.pictex}}\close + \frac{3}{5}\open{\input{e1e.pictex}}\close + \frac{4}{15}\open{\input{b1.pictex}}\close. 
\end{split}
\end{equation}
\end{lemma}

Five new strata with smooth genus 2 components appear in this relation:
$$ \open{\input{d.pictex}}\close, \open{\input{d1.pictex}}\close, \open{\input{d2.pictex}}\close, \open{\input{g.pictex}}\close, \open{\input{capee.pictex}}\close. $$
A Weierstrass point condition is denoted by a $W$ on the node or
marking (and is always on the genus 2 component).
The letters $x, \barr{x}$ designate a hyperelliptic conjugate pair
condition --- they are {\em not} marking labels. 

On the new strata with elliptic components,
$$ \open{\input{b4.pictex}}\close, \open{\input{b3.pictex}}\close, \open{\input{c5.pictex}}\close, $$
the letters $x, y$ designate an imposed linear equivalence on
the markings and nodes: the divisor sum of the points lettered $x$
must be equivalent to the divisor sum of the points lettered $y$.
In the first two strata, the sum of the points lettered
$x$ must be equivalent to twice the node point on the left component. 
The third stratum has an imposed linear equivalence on the
normalization: the sum of the two preimages of the node must be
equivalent to the sum of the points lettered $x$. Our previous
convention regarding unassigned markings holds: the sum over the 3!
marking choices is taken.

Theorem 1 will be obtained from Lemma \ref{lemma1} by expressing all
occurring strata in terms of descendent strata.

\begin{proof}[Proof of Lemma \ref{lemma1}]
The method of proof is by direct calculation of the map
$\lambda_* \pi^*$ on each term of (\ref{wdvv1}). The technique
is identical to the elliptic calculations in \cite{rahul}.
We will give a representative example of the computation. 

Let $C$ be the codimension 2 stratum of $\Mbar_{0,3+6}$ occurring as
the divisor $D(\star p_3|b_1\etc b_6)$ in $\Mbar_{0,A}$. The
generic point of $C$ corresponds to a genus 0 stable curve which is a
chain of three rational components $U,V,W$. The
marking distribution is as follows:
$p_1$ and $p_2$ on $U$, $p_3$ on $V$, and 
the six branch points $b_1,\etc,b_6$ on $W$. An admissible
cover over such a curve consists of disjoint \'etale double 
covers of $U$ and $V$ and a genus 2 double cover of $W$ branched over 
$b_1,\etc,b_6$.

\vspace{5pt}
\centerline{\input{admiss.pictex}}
\vspace{2pt}

The preimage $\pi^{-1}(C)$ has four irreducible components
$\Sigma_1,\etc,\Sigma_4$ corresponding to four ways (up to isomorphism
of the cover) of distributing the marked points among the sheets of 
the cover. The component $\Sigma_4$ parametrizes covers with all three
marked points placed on the same sheet, whereas each of the components
$\Sigma_i,\; 1\leq i\leq 3$ 
parametrizes covers with $p_i$ placed on one sheet and the remaining
two markings placed on the other sheet. Since the projection $\pi$ is a
finite group quotient, we can compute the pull-back of the class of
$C$ by $\pi^*$ using the following lemma (see \cite{vistoli} for
the proof).

\begin{lemma}
\label{pbck}
Let $\G$ be a finite group, let $X$ be an irreducible algebraic
variety with a $\G$-action, and let $\alpha$ denote the quotient
morphism $\alpha: X \to X/\G$. There exists a pull-back $\alpha^*:
A_*(X/\G) \to A_*(X)$ defined by: 
$$\alpha^*[V]=  {|\text{\em Stab}(V)|}\cdot
{[\alpha^{-1}(V)_{\red}]},$$
where $V$ is an irreducible subvariety of $X/\G$, the scheme
$\alpha^{-1}(V)_{\red}$ is the reduced preimage of $V$, and
$|\text{\em Stab}(V)|$ is the size of the generic stabilizer of points
over $V$.  
\end{lemma}

The stabilizer of points over $C$ is trivial. Hence, by Lemma
\ref{pbck}, we get 
$$\pi^*[C] = \Sigma_1+\Sigma_2+\Sigma_3+\Sigma_4.$$
The $\lambda$ push-forward of the class of the component $\Sigma_4$ is
$6!\times \open{\input{tr1ord.pictex}}\close$, whereas the push-forwards of the classes of
$\Sigma_1, \Sigma_2, \Sigma_3$ are
$$6!\times\open{\input{capeeord1.pictex}}\close,\;\; 6!\times\open{\input{capeeord2.pictex}}\close\;\text{and}\;\;
6!\times\open{\input{capeeord3.pictex}}\close$$
respectively.
Let $C_\sym$ denote the cycle in $A^2(\Mbar_{0,3+6})$ obtained by
symmetrizing $C$ with respect to the $\sy_3$-action: $C_\sym =
\sum_{\sigma\in\sy_3} \sigma^*(C)$. We obtain 
$$\lambda_*\pi^* [C_\sym] = 6!\times\open{\input{tr1.pictex}}\close + 3\cdot 6!\times\open{\input{capee.pictex}}\close.$$
The push-forwards of the other cycles are computed similarly.
\end{proof}

\subsection{The second equation}
\label{firstred}
We find an expression for the codimension 2 class $\open{\input{g.pictex}}\close$ 
appearing in (\ref{rel1}). 
Consider the boundary divisor $D=D(p_3b_6|p_1p_2b_1\ldots b_5)$
in $\Mbar_{0,3+6}$. It is isomorphic to $\Mbar_{0,A}$ with the marking
set $A=\{\star,p_1,p_2,b_1,\etc, b_5\}$.
Consider the following 4-point linear equivalence on $\Mbar_{0,A}$: 
$$ (p_1p_2|b_1\star)\sim (p_1b_1|p_2\star).$$ 
As in Section \ref{firsteq}, let $\Lambda_D$ be the corresponding
relation in $A^2(\Mbar_{0,3+6})$.

\begin{lemma}
\label{lemfirst}
The application of $\lambda_*\pi^*$ to relation $\Lambda_D$ yields
(\/$2\cdot 5!$ times) the following rational equivalence in
$A^2(\Mbar_{2,3})$: 
\begin{equation}
\label{rel2}
\begin{split}
&\open{\input{gord.pictex}}\close = \frac{2}{5}\open{\input{bubbaord.pictex}}\close + 2\open{\input{d2ord.pictex}}\close - \open{\input{d1ord.pictex}}\close \\[3pt]
&+ \frac{6}{5}\open{\input{a3ord.pictex}}\close - \frac{4}{5}\open{\input{a2ord.pictex}}\close + \frac{4}{5}\open{\input{b1ord.pictex}}\close -
\frac{1}{5}\open{\input{b2ord.pictex}}\close.  
\end{split}
\end{equation}
\end{lemma}

\subsection{The third equation}
We find an expression for the codimension 2 class $\open{\input{bubbaord.pictex}}\close$ 
appearing in (\ref{rel2}). 
Consider the boundary divisor $D=D(p_3b_6|p_1p_2b_1\ldots b_5)$ 
in $\Mbar_{0,3+6}$. It is isomorphic to $\Mbar_{0,A}$ with the marking
set $A=\{\star,p_1,p_2,b_1,\etc,b_5\}$. 
Consider the following 4-point linear equivalence on $\Mbar_{0,A}$: 
$$(p_1b_1|\star b_2) \sim ( p_1\star|b_1b_2).$$
Let $\Lambda_D$ be the corresponding relation in $A^2(\Mbar_{0,3+6})$. 

\begin{lemma}
The application of $\lambda_*\pi^*$ to relation $\Lambda_D$ yields
(\/$4\cdot 4!$ times) the following rational equivalence in
$A^2(\Mbar_{2,3})$: 
\begin{equation}
\label{rel3}
\begin{split}
&\open{\input{bubbaord.pictex}}\close = 5\open{\input{d2ordd.pictex}}\close + 5\open{\input{dord.pictex}}\close + \frac{3}{2}\open{\input{a4ord.pictex}}\close +
\frac{3}{2}\open{\input{a3ordd.pictex}}\close - \frac{3}{2}\open{\input{a3ord.pictex}}\close  \\[3pt]
&- \frac{3}{2}\open{\input{a2ord.pictex}}\close + \frac{3}{2}\open{\input{b0ord.pictex}}\close + \frac{3}{2}\open{\input{b1ordd.pictex}}\close -
\frac{1}{2}\open{\input{b1ord.pictex}}\close - \frac{1}{2}\open{\input{b2ord.pictex}}\close + \frac{1}{8}\open{\input{b5ord.pictex}}\close.   
\end{split}
\end{equation}
\end{lemma}

\subsection{Weierstrass to descendent stratum classes}
Combining relations (\ref{rel1}), (\ref{rel2}) and (\ref{rel3}) yields
a relation in $A^2(\Mbar_{2,3})$ in which the only strata with a
smooth genus 2 component of the corresponding generic stable curve
are:  
\begin{equation}
\label{weier}
\open{\input{tr1.pictex}}\close, \open{\input{d.pictex}}\close, \open{\input{d1.pictex}}\close, \open{\input{d2.pictex}}\close, \open{\input{capee.pictex}}\close. 
\end{equation}
The first stratum is pure boundary. Expressions for the next 3 classes
in terms of descendent stratum classes are obtained using the
following result.   

\begin{lemma}
The following linear equivalence holds in $A^1(\Mbar_{2,1})$: 
\begin{equation*}
\open{\input{wei.pictex}}\close = 3\psi_1 -\frac{6}{5}\open{\input{del1.pictex}}\close - \frac{1}{10}\open{\input{del0.pictex}}\close.
\end{equation*}
\end{lemma}

A proof can be found in \cite{eizhar}.
The relations obtained are:
\begin{align*}
\open{\input{d.pictex}}\close   &= 3\open{\input{p.pictex}}\close - \frac{6}{5}\open{\input{ee.pictex}}\close -\frac{1}{10}\open{\input{en.pictex}}\close, \\
\open{\input{d1.pictex}}\close &= 3\open{\input{p1.pictex}}\close - 3\open{\input{tr1.pictex}}\close - \frac{6}{5}\open{\input{e1e.pictex}}\close - \frac{6}{5}\open{\input{ee1.pictex}}\close -
\frac{1}{10}\open{\input{en1.pictex}}\close, \\
\open{\input{d2.pictex}}\close &= 3\open{\input{p2.pictex}}\close - 3\open{\input{tr1.pictex}}\close - \frac{6}{5}\open{\input{e1e.pictex}}\close - \frac{6}{5}\open{\input{ee1.pictex}}\close -
\frac{1}{10}\open{\input{en1.pictex}}\close.
\end{align*}

Expression for the last class in (\ref{weier}) in terms of descendent
stratum classes is obtained by combining the above formulas with the
following result.

\begin{lemma}
\label{lemlast}
The following linear equivalence holds in $A^1(\Mbar_{2,2})$:
\begin{equation}
\label{conjugate}
\open{\input{conj.pictex}}\close = \frac{2}{3}\open{\input{w2.pictex}}\close  - \open{\input{last2.pictex}}\close +
\frac{3}{5}\open{\input{last1.pictex}}\close - \frac{2}{5}\open{\input{last12.pictex}}\close -
\frac{1}{30}\open{\input{last0.pictex}}\close. 
\end{equation}
\end{lemma}
\begin{proof}
Consider the following 4-point linear equivalence on $\Mbar_{0,2+6}$:  
\begin{equation}
\label{twotwo}
(p_1p_2|b_1b_2) \sim ( p_1b_1|p_2b_2). 
\end{equation}

Consider the morphisms
$$\begin{CD}
\Hbar_{2,2,2} @>{\lambda}>> \Mbar_{2,2} \\ 
@V{\pi}VV \\
\Mbar_{0,2+6} \\ 
\end{CD}$$
The application of $\lambda_*\pi^*$ to relation (\ref{twotwo}) yields 
$3\cdot 5!$ times relation (\ref{conjugate}).
\end{proof}

\subsection{Genus 1}
We list the relations in $A^1(\Mbar_{1,n})$ necessary to
express the classes in (\ref{rel1}), (\ref{rel2}) and (\ref{rel3})
with genus 1 components constrained by linear equivalence conditions
in terms of the boundary classes. These relations are obtained
similarly to the above relations, only using the space of elliptic 
admissible double covers (with 4 branch points). See \cite{pasha} for
details. 

\begin{align*}
\open{\input{g12.pictex}}\close &= 3\open{\input{two1.pictex}}\close + \frac{1}{4}\open{\input{two0.pictex}}\close, \\
\open{\input{g13.pictex}}\close &= 2\open{\input{three4.pictex}}\close + 2\open{\input{three2.pictex}}\close + 2\open{\input{three3.pictex}}\close - \open{\input{three1.pictex}}\close +
\frac{1}{6}\open{\input{three0.pictex}}\close, \\
\open{\input{g14.pictex}}\close &= \open{\input{fourall.pictex}}\close +\frac{1}{6}\open{\input{four1un.pictex}}\close + \open{\input{four13.pictex}}\close + \open{\input{four14.pictex}}\close 
     + \open{\input{four24.pictex}}\close + \open{\input{four23.pictex}}\close \\
    &- \open{\input{four12.pictex}}\close - \open{\input{four34.pictex}}\close +  \frac{1}{12}\open{\input{four0.pictex}}\close.
\end{align*}

Lemmas \ref{lemfirst}--\ref{lemlast} and the above equations express
all the strata occurring in relation (\ref{rel1}) in terms of
descendent stratum classes. The resulting relation among descendent
stratum classes is (\ref{maine}). The proof of Theorem 1 is complete.

\section{\bf Descendent integrals on $\P^2$}

\subsection{1-cotangent line integrals}
\label{1cot}
Let $X$ be a nonsingular projective variety. Define the 
{\em 1-cotangent line descendents} to be the integral invariants
(\ref{grav}) of $X$ with at most 1 cotangent line class. 
By Getzler's topological recursion relations in genus 2
(\cite{getzler2}), the 1-cotangent line integrals determine
all genus 2 descendent integrals. 
Relation (\ref{maine}) yields new universal equations satisfied by
descendents invariants of $X$.
It was initially hoped these equations  
would show that genus 2 descendents could
be uniquely reconstructed from Gromov-Witten
invariants for any space $X$. 
Such universal reconstruction results are known in genus
0 and 1. 
However, we found all of our reconstruction strategies thwarted
by specific unexpected linear relations among the
coefficients of relation (\ref{maine}).
While universal descendent reconstruction may hold in 
genus 2, new ideas are required: either
subtle strategies involving relation (\ref{maine}) or
yet another genus 2 relation.

Much more can be said if attention is restricted to specific target
spaces. In this section, relation (\ref{maine}) is shown to determine  
all 1-cotangent line descendents of $\P^2$ from degree 0 ones
and the lower genus invariants. Hence, using  Getzler's topological 
recursion relations, all genus 2 descendent invariants of $\P^2$ are
determined via descendent stratum relations in $\Mbar_{2,n}$. 

Let $T_0, T_1, T_2$ be the standard cohomology basis of
$\P^2$ given by the fundamental, hyperplane, and point classes
respectively. 
On $\P^2$, the basic 
1-cotangent line descendent integrals assemble in 3 series:
\begin{equation}
\label{threeser}
\begin{split}
&\N_d  = \lan T_2^{3d+1} \ran _{2,d}, \\
&\HH_d = \lan (T_1\psi) \cdot T_2^{3d} \ran_{2,d}, \\
&\PP_d = \lan (T_2\psi) \cdot T_2^{3d-1}\ran_{2,d}. 
\end{split}
\end{equation}
All descendent integrals with at most 1 cotangent
line class may be reduced to one of the basic integrals
above via the string, dilaton, and divisor equations.
The first two series are defined in degree 0.
A direct virtual class computation yields:
\begin{equation}
\label{valss}
\begin{split}
\N_0  &= \lan T_2 \ran_{2,0} =0, \\
\HH_0 &= \lan T_1\psi \ran_{2,0}= 
-6\cdot \int_{\Mbar_{2}} \lambda_1 \lambda_2 = - \frac{1}{960}.
\end{split}
\end{equation}
The integral of the Chern classes of the Hodge bundle over
$\Mbar_2$ may be computed by a method due to Faber (\cite{faber}). 

Following conventions, let the variable $t^i_j$ correspond to the
class $T_i \psi^j$, where $0\leq i \leq 2$ and $j \geq 0$. 
Let $t$ denote the variable set $\{t^i_j\}$, and let
$\gamma= \sum t^i_j T_i \psi^j$ be the formal sum.
Let $F_{2, \proj^2} (t)$ denote the full genus 2 gravitational
potential function: 
\begin{equation}
\label{bigf}
F_{2, \proj^2} (t)= 
\sum_{d\geq 0} \sum_{n\geq 0} \frac{1}{n!} \lan \gamma^n \ran_{2,d}.
\end{equation}
We will consider the cut-off of $F_{2,\proj^2}(t)$ at 1-cotangent
line: 
$$ 
F^1_{2, \proj^2}(t^0_0, t^1_0, t^2_0, t^0_1, t^1_1, t^2_1) =
F_{2, \proj^2}(t)|_{\{t^i_j=0 \ \forall j\geq 2,\ t^i_1t^j_1=0 \
  \forall i,j\}}. 
$$
The function $F_{2, \P^2}^1$ may be written explicitly in terms
of the 3 basic series:
\begin{eqnarray*}
F^1_{2, \proj^2}& = & 
\ \ \ \sum_{d\geq 0} \N_d e^{dt^1_0}
\frac{(t^2_0)^{3d+1}}{(3d+1)!} \\
& & + \sum_{d\geq 0} \N_d d\, t^0_0 t^1_1 e^{dt^1_0}
\frac{(t^2_0)^{3d+1}}{(3d+1)!} 
+ \sum_{d\geq 0} \N_d t^0_0 t^2_1 e^{dt^1_0}
\frac{(t^2_0)^{3d}}{(3d)!} \\ 
& & 
+ \sum_{d\geq 0} \N_d (3d+3+d\, t^1_0)\ t^0_1\  e^{dt^1_0}
\frac{(t^2_0)^{3d+1}}{(3d+1)!} \\  
& & + \sum_{d\geq 0} (\HH_d+ \N_dt^1_0)\ t^1_1 \ e^{dt^1_0} 
\frac{(t^2_0)^{3d}}{(3d)!}\\
& & + \sum_{d\geq 1} \PP_d\ t^2_1\ e^{dt^1_0} \frac{(t^2_0)^{3d-1}}{(3d-1)!}.
\end{eqnarray*}
Relation (\ref{maine}) will yield differential equations which determine
$F^1_{2, \proj^2}$ from degree 0 values (\ref{valss}) and the
known elliptic and rational Gromov-Witten potentials.

The descendent integrals (\ref{grav}) are defined via the cotangent
lines on the moduli space of maps. Let $2g-2+n >0$, and let 
$$ f: \Mbar_{g,n}(X,\beta) \to \Mbar_{g,n} $$
be the forgetful map. The following integrals are similar to the
descendents: 
\begin{equation}
\label{downd} 
\underset{[\Mbar_{g,n}(X, \beta)]^\virt}{\int}  
 \ev_1^*(\gamma_1) \scup f^*(\psi_1^{\alpha_1})\scup \cdots \scup
 \ev_n^*(\gamma_n) \scup f^*(\psi_n^{\alpha_n}). 
\end{equation}
The cotangent lines here are pulled back from the
moduli space of stable pointed curves.
The integral invariants (\ref{downd}) are naturally equivalent
to the descendent invariants. The two sets of invariants are related
by universal invertible transformations (see \cite{kontsmanin2}). 
The differential equations obtained from relation (\ref{maine})
via the splitting formula in Gromov-Witten theory are
most conveniently expressible in terms of generating functions for
the invariants (\ref{downd}). However, in the
restricted case of 1-cotangent line integrals on $\P^2$, the two sets
of invariants are exactly equal.

\begin{lemma}
\label{lhap}
Let $g>0$ and  $d\geq 0$. 
Let $\gamma_i \in H^*(\P^2, \Q)$.
An equality of 1-cotangent line integrals holds for $\P^2$:
\begin{equation}
\label{happy} 
\underset{[\Mbar_{g,n}(\P^2, d)]^\virt}{\int} f^*(\psi_1) \scup
\prod_{i=1}^n \ev_i^*(\gamma_i)\ =
\underset{[\Mbar_{g,n}(\P^2, d)]^\virt}{\int} \psi_1 \scup
\prod_{i=1}^n \ev_i^*(\gamma_i).
\end{equation}
\end{lemma}

\begin{proof}
Define the pairing matrix $(g_{ij})$ by 
$$ g_{ij}=\int_{\P^2} T_i\scup T_j. $$
Let $(g^{ij})=(g_{ij})^{-1}$, so that $\sum g^{ij} T_i \boxtimes T_j$
is the class of the diagonal in $\P^2 \times \P^2$. 
By the formulas of \cite{kontsmanin2}, the difference between the
integrals (\ref{happy}) is expressed in terms of pure
Gromov-Witten invariants:
$$ \sum \lan \gamma_1 \cdot T_i\ran_{0,d_1} \  g^{ij}\  \lan T_j 
\cdot \gamma_2 \cdots \gamma_n\ran_{g, d_2}, $$
where the sum is over all degree splittings $d_1+d_2=d$  satisfying
$d_1>0$ and the diagonal splitting. The only nonvanishing 2-point
genus 0 Gromov-Witten invariant of positive degree is $\lan T_2 \cdot
T_2 \ran_{0,1}$. Hence, only the summands with $j=0$ can
contribute. However, Gromov-Witten invariants of genus $g>0$ with an 
argument $T_0$ vanish by the axiom of the fundamental class.   
\end{proof}

Let $\widetilde{F}_{2, \proj^2}$ denote the full genus 2 potential
function defined using the integrals (\ref{downd}) instead of the
descendents, and let $\widetilde{F}^1_{2, \proj^2}$
be the 1-cotangent line cut-off. Then, by Lemma \ref{lhap},
$F^1_{2,\proj^2}= \widetilde{F}^1_{2, \proj^2}$.

\subsection{Pull-backs of descendent stratum classes}
\label{pullbacks}
In order to find differential equations via the splitting formula, 
the pull-backs of
relation (\ref{maine}) to $\Mbar_{2,3+l}$ will be required. 
More generally, let $\nu$ denote the forgetful map,
\begin{equation}
\label{fgful}
\nu:\Mbar_{g,n+l} \to \Mbar_{g,n}.
\end{equation}
The pure boundary strata (without cotangent line classes) of
$\Mbar_{g,n}$ are naturally indexed by stable, $n$-pointed, genus $g$
dual graphs: the vertices, edges, and 
markings of the graph correspond to the components, nodes, and marked
points of the curve with moduli point generic in the stratum. The dual
graph of a stratum is equivalent data to the topological type and
marking distribution of its generic curve (the latter is the notation
chosen in Section \ref{relcomp} for visual simplicity).
Let $\Gamma$ be a stable dual graph. 
The class of the stratum $S_\Gamma \subset \Mbar_{g,n}$
pulls back under the forgetful map $\nu$ 
to the sum of classes of boundary strata obtained by all possible
distributions of the $l$ extra points $\{q_1,\etc,q_l\}$ on the
vertices of $\Gamma$. 
For example, the pull back of the class $\open{\input{tr1ord.pictex}}\close$ from $\Mbar_{2,3}$
to $\Mbar_{2,3+l}$ is given by:
\begin{equation}
\label{pbacky}
\nu^* \open{\input{tr1ord.pictex}}\close = \sum_{A\disj B\disj C = \{q_1,\etc,q_l\}} \Bigg[\raisebox{-4ex}
{\input{pbktr1.pictex}}\hskip-4pt\Bigg].
\end{equation}
While the pull-back of a descendent stratum class is a sum
of descendent stratum classes, the pull-back is not obtained
simply by distributing the extra markings as in (\ref{pbacky}). 
The additional complexity occurs because the pull-back of
a cotangent line class via the forgetful map (\ref{fgful}) is 
{\em not} the corresponding cotangent line class on the domain. 
\begin{lemma}
\label{discrep}
There is an equality:
$$ 
\nu^*(\psi_1) = \psi_1 - \sum_{A\disj B=\{q_1,\etc,q_l\}}\Bigg[\hskip2pt\raisebox{-3.5ex}
  {\input{psipull.pictex}}\hskip-2pt\Bigg], 
$$
where the sum is over all \underline{stable} distributions of 
the $l$ extra points.
\end{lemma}
Lemma \ref{discrep} easily implies the pull-back formulas for
the three codimension 2 descendent strata in $\Mbar_{2,3}$:
\begin{equation}
\label{pbks}
\begin{split}
\nu^*\open{\input{pord.pictex}}\close &= \sum_{A\disj B = \{q_1,\etc,q_l\}} \Bigg[\raisebox{-3.5ex}
  {\input{pbkp.pictex}}\hskip-4pt\Bigg]  
- \sum_{\substack{A\disj B\disj C = \{q_1,\etc,q_l\} \\ B\neq \emptyset}}
\Bigg[\raisebox{-4ex}
{\input{pbkpcorr.pictex}}\hskip-4pt\Bigg], \\ 
\nu^*\open{\input{p1ord.pictex}}\close &= \sum_{A\disj B = \{q_1,\etc,q_l\}} \Bigg[\raisebox{-3.5ex}
  {\input{pbkp1.pictex}}\hskip-4pt\Bigg]  
- \sum_{\substack{A\disj B\disj C = \{q_1,\etc,q_l\} \\ B\neq \emptyset}}
\Bigg[\raisebox{-4ex}
{\input{pbkp1corr.pictex}}\hskip-4pt\Bigg], \\
\nu^*\open{\input{p2ord.pictex}}\close &= \sum_{A\disj B = \{q_1,\etc,q_l\}} \Bigg[\raisebox{-3.5ex}
  {\input{pbkp2.pictex}}\hskip-4pt\Bigg]  
- \sum_{\substack{A\disj B\disj C = \{q_1,\etc,q_l\} \\ B\neq \emptyset}}
\Bigg[\raisebox{-4ex}
{\input{pbkp2corr.pictex}}\hskip-4pt\Bigg]. 
\end{split}
\end{equation}

Using these formulas, the pull-backs of relation (\ref{maine}) to
$\Mbar_{2,3+l}$ yield descendent stratum class relations.

\subsection{The splitting formula in Gromov-Witten theory}
We review here the splitting formula following \cite{kontsmanin1}.
Let $X$ be a nonsingular projective variety with cohomology basis
$T_0, \etc, T_m$. Let $g_{ij}=\int T_i\scup T_j$, and let
$(g^{ij})=(g_{ij})^{-1}$.
Let $2g-2+n>0$. As in Section (\ref{1cot}), let
$$ f: \Mbar_{g,n}(X,\beta) \to \Mbar_{g,n} $$
be the forgetful map. 
Let $\gamma_i \in H^*(X, \Q)$.
A {\em Gromov-Witten class} $\lan \prod_{i=1}^{n} \gamma_i
\ran_{g,\beta}$ is an element of $H^*(\Mbar_{g,n}, \Q)$ defined via
push-forward by $f$:
$$ \lan\gamma_1 \cdots \gamma_n\ran_{g,\beta} = 
f_* \Big( [\Mbar_{g,n} (X, \beta)]^\virt \scap \prod_{i=1}^n
\ev_i^*(\gamma_i) \Big). $$

Let $\Gamma$ be a stable, $n$-pointed, genus $g$ dual graph
corresponding to a stratum class $[S_\Gamma]$. Let $\mathcal{V}$ and
$\mathcal{E}$ denote the set of vertices and edges of $\Gamma$. 
Let $\mathcal{E}_{\flag}$ denote the set of edge flags
(an edge is made up of 2 edge flags).
For each vertex $v\in \mathcal{V}$,
let $\mathcal{I}(v)$ 
and $\mathcal{E}_{\flag}(v)$ 
be the set of markings and edge flags incident at $v$.
Let $\text{val}(v)$ denote the valence of $v$ in $\Gamma$, and let
$g(v)$ be the genus assignment. For an edge $E\in \mathcal{E}$, let 
$E_f$ and $E_f^*$ denote the corresponding two edge flags.
Let $\text{A}$ be the automorphism group of the
graph (as a pointed graph with genus assignment). Let $\iota: S_\Gamma 
\into \Mbar_{g,n}$ denote the inclusion. Let
$$\pi: \prod_{v\in V} \Mbar_{g(v), \text{val}(v)}
\to S_\Gamma$$
be the natural map.
The splitting formula is:
$$
\iota^* \lan \prod_{i=1}^{n} \gamma_i\ran_{g,\beta} =
\frac{1}{|\text{A}|}\  \pi_* \Big(
\sum_{\epsilon, \phi }
\prod_{E\in \mathcal{E}} g^{\epsilon(E_f) \epsilon(E_f^*)} 
\prod_{v\in \mathcal{V}}
\ \ \lan \prod_{i\in \mathcal{I}(v)} \gamma_i \prod_{e\in
  \mathcal{E}_{\flag}(v)} T_{\epsilon(e)}\ran_{g(v), \phi(v)} 
\Big). 
$$
The sum is over all functions
 $\epsilon: \mathcal{E}_{\flag} \to
\{0, \etc, m\}$ from edge flags to the basis index set, and
all functions $\phi:\mathcal{V} \to H_2(X, \Z)$
satisfying  $\sum_{v\in \mathcal{V}} \phi(v) = \beta$.

\subsection{Differential equations}
Relation (\ref{maine}), the pull-back formulas of Section
\ref{pullbacks}, and the splitting formula naturally yield
differential equations for the potential function
$\widetilde{F}_{2,X}$. Let $\widetilde{F}_{0,X}$ and
$\widetilde{F}_{1,X}$ denote the full potential functions in genus 0
and 1 with respect to the integrals (\ref{downd}). Let
$\tilde{\gamma}= \sum t^i_j T_i f^*(\psi^j)$ be the formal sum. Then
\begin{eqnarray*}
\widetilde{F}_{0,X}(t) & = &  \sum_{d\geq 0} \sum_{n\geq 3} \frac{1}{n!}
\lan \tilde{\gamma}^n 
\ran_{0,d}, \\
\widetilde{F}_{1,X}(t)& = &  \sum_{d\geq 0} \sum_{n\geq 1} \frac{1}{n!}
\lan \tilde{\gamma}^n \ran
_{1,d}.
\end{eqnarray*}
These generating functions are sums over the stable range $\{2g-2+n >0\}$
of {$n$-pointed curves of genus $g$. Let $\widetilde{F}_{0,X}^0$ and
$\widetilde{F}_{1,X}^0$ denote the restriction to the small phase
space $\{t^i_j=0 \ \forall j\geq 1\}.$
These are the 0-cotangent line cut-offs and involve only the
Gromov-Witten invariants. 

The differential equations are now described. An  equation is obtained
for every assignment of variables $\{t^i_j\}$ to the 3 markings. 
Fix such an assignment 
$(t^{i_1}_{j_1}, t^{i_2}_{j_2}, t^{i_3}_{j_3})$.
A differential equation
\begin{equation}
\label{difff}
\mathcal{D}^{i_1i_2i_3}_{j_1j_2j_3} (\widetilde{F}_{0,X},
\widetilde{F}_{1,X}, \widetilde{F}_{2,X})=0
\end{equation}
is constructed from relation (\ref{maine}) in the following manner.

A pure boundary stratum $S_\Gamma \subset \Mbar_{2,3}$ naturally
yields a differential expression: place potentials on the vertices of
the dual graph, insert the 3 point conditions via differentiation,
sum over diagonal splittings at the edges, and divide by the number of
graph automorphisms.  For example, the  strata $\open{\input{tr1ord.pictex}}\close$ and
$\open{\input{er3ord.pictex}}\close$ yield the respective expressions
\begin{gather*}
  \frac{\partial \widetilde{F}_{2,X}}{\partial t^e_0} g^{ef}
  \frac{\partial^3 \widetilde{F}_{0,X}}{\partial t^f_0 \partial 
  t^{i_3}_{j_3} \partial t^k_0} g^{kl} \frac{\partial^3
  \widetilde{F}_{0,X}}{\partial t^l_0 \partial t^{i_1}_{j_1} \partial
  t^{i_2}_{j_2}}, \\
\frac{1}{2}\cdot \frac{\partial \widetilde{F}_{1,X}}{\partial
  t^e_0} g^{ef} \frac{\partial^5 \widetilde{F}_{0,X}}{\partial t^f_0
  \partial t^{i_1}_{j_1} \partial t^{i_2}_{j_2} \partial t^{i_3}_{j_3}
  \partial t^k_0} g^{kl} \frac{\partial\widetilde{F}_{1,X}}{\partial
  t^l_0}.  
\end{gather*}

Each descendent stratum in relation (\ref{maine}) yields a two-term
differential expression. The first term is again obtained by placing
potentials on the vertices, inserting point conditions, and summing
over diagonal splittings (no automorphisms occur for these
graphs). For example, the first terms for the descendent strata
$\open{\input{p2ord.pictex}}\close, \open{\input{p1ord.pictex}}\close$ are
$$ 
\frac{\partial^2 \widetilde{F}_{2,X}}{\partial t^{i_1}_{j_1} 
\partial t^e_1} g^{ef} 
\frac{\partial^3 \widetilde{F}_{0,X}}{\partial t^f_0 \partial
t^{i_2}_{j_2} \partial t^{i_3}_{j_3}}\ \ \text{and}\ \ 
\frac{\partial^2 \widetilde{F}_{2,X}}{\partial t^{i_1}_{j_1+1} 
\partial t^e_0} g^{ef} 
\frac{\partial^3 \widetilde{F}_{0,X}}{\partial t^f_0 \partial
t^{i_2}_{j_2} \partial t^{i_3}_{j_3}}. 
$$
The second term is obtained from the correction graphs in the
pull-back formulas (\ref{pbks}). For $\open{\input{p2ord.pictex}}\close$, the correction 
graph is $\bigg(\hskip2pt\raisebox{-2ex}
{\input{corrgraph.pictex}}\hskip-2pt\bigg)$. The second term it yields is
$$
- \frac{\partial^2 \widetilde{F}_{2,X}}{\partial t^{i_1}_{j_1} 
\partial t^e_0} g^{ef}
 \frac{\partial^2 \widetilde{F}_{0,X}}{\partial t^f_0 
\partial t^k_0} g^{kl}
\frac{\partial^3 \widetilde{F}_{0,X}}{\partial t^l_0 \partial t^{i_2}_{j_2}
\partial t^{i_3}_{j_3}}.
$$

Equation (\ref{difff}) is constructed by replacing the strata classes
in (\ref{maine}) by the corresponding differential expressions.
Equation (\ref{difff}) is then easily proven by pulling-back
(\ref{maine}) to the moduli spaces and applying the splitting formula. 

Equations strictly among 1-cotangent line integrals of $X$ may be
obtained from the differential equations (\ref{difff}) by the
following method. Let the variables $t^i_j$ assigned to the 3 markings 
correspond to pure cohomology classes (i.e. $j=0$).
Restrict the left side of (\ref{difff}) to the small phase space:
\begin{equation}
\label{smdifff}
\mathcal{D}^{i_1i_2i_3}_{0\; 0\; 0} (\widetilde{F}_{0,X},
\widetilde{F}_{1,X}, \widetilde{F}_{2,X}) \ |_{\{t^i_j=0 \ \forall
  j\geq 1\}}=0.  
\end{equation}
The derivatives $\partial/\partial t^i_j$ with $j\geq 1$
occur only in the terms of $\mathcal{D}$ obtained from the strata: 
$$
\open{\input{p.pictex}}\close, \open{\input{p1.pictex}}\close, \open{\input{p2.pictex}}\close.
$$
In each of these terms, the derivative appears only once, with
$j=1$. Moreover, the derivative acts on the genus 2 potential. 
Hence, equation (\ref{smdifff}) implies:
\begin{equation}
\label{ssmdifff}
\mathcal{D}^{i_1 i_2 i_3}_{0\; 0\; 0}
(\widetilde{F}^0_{0,X},
\widetilde{F}^0_{1,X}, \widetilde{F}^1_{2,X})    \ |_{\{t^i_j=0 \
  \forall j\geq 1\}}=0. 
\end{equation}
The coefficient relations obtained from (\ref{ssmdifff}) are exactly
among 1-cotangent line integrals of genus 2 and Gromov-Witten
invariants of genus 0 and 1.

\subsection{Recursions for $\P^2$}
The differential equations (\ref{ssmdifff}) yield recursive
relations among the 3 basic series (\ref{threeser}) of 1-cotangent
line integrals of $\P^2$ (involving the lower genus
Gromov-Witten invariants). Each choice of
point assignment provides such recursions.
In each degree $d$, there are 3 basic 1-cotangent line
integrals. Hence, 3 independent recursions are
required. Four distinct equations of the form (\ref{ssmdifff}) are
obtained by the marking choices:
\begin{gather*}
(T_1,T_1,T_1), \\
(T_1,T_1,T_2), \\
(T_1,T_2,T_2), \\
(T_2,T_2,T_2). 
\end{gather*}
The relations are computed and are easily seen to determine
the 3 series from the degree 0 values (\ref{valss}) and
Gromov-Witten invariants of genus 0 and 1.

We include here the recursion obtained from the assignment
$(T_1,T_1,T_1)$:
\begin{equation*}
\begin{split}
-3 & \HH_d + 3d\PP_d
= \sum_{\substack{d_1+d_2+d_3 = d \\ d_i>0}}
\Big(p_{200}\;\N_{d_1}\R_{d_2}\R_{d_3}
+p_{110}\;\E_{d_1}\E_{d_2}\R_{d_3}\Big) \\
&+\sum_{\substack{d_1+d_2=d \\ d_i>0}} \Big(p_{20}\;\N_{d_1}\R_{d_2}
+p_{h0}\;\HH_{d_1}\R_{d_2} + p_{11}\;\E_{d_1}\E_{d_2}
+p_{10}\;\E_{d_1}\R_{d_2}\Big) \\
&{\textstyle
-\frac{1}{960}d^4(d-1)(d-2)\R_d+\frac{1}{40}d^2(5d-6)\E_d}.
\end{split}
\end{equation*}
\vspace{2pt}

The polynomial coefficients are defined by
\begin{align*}
&p_{200} ={\textstyle -2\binom{3d-1}{3d_1+1,\;3d_2-1,\;3d_3-1}
d_1d_2^2d_3^3(d_2+d_3)},\\ 
&p_{110} ={\textstyle \frac{1}{5}\binom{3d-1}{3d_1,\;3d_2,\;3d_3-1} 
d_1d_2d_3^3 (-9d_1d_2+6d_2^2-12d_2d_3+d_3^2)}, \\
&p_{20} ={\textstyle
\binom{3d-1}{3d_1+1} d_2(3d_1^2-10d_1d_2+4d_2^2)+
\binom{3d-1}{3d_1} d_2^2(3d_1+2d_2)},\\ 
&p_{h0} ={\textstyle 2\binom{3d-1}{3d_1} d_2^4},\\
&p_{11} ={\textstyle \frac{3}{5}\binom{3d-1}{3d_1}
  d_1(4d_1^2-9d_1d_2+2d_2^2)},\\
&p_{10} ={\textstyle \frac{1}{120}\binom{3d-1}{3d_1} 
  d_1d_2^3 (-18d_1^2+36d_1d_2-6d_2^2+5d_1^3-33d_1^2d_2
  +3d_1d_2^2+d_2^3)}.
\end{align*}
 
A calculation of the first 10 values of the 3 series is tabulated
below. There are at least four other mathematical methods to obtain
the series $\N_{d}$: the degenerations of \cite{ran} and
\cite{caphar}, the hyperelliptic methods of \cite{graber}, and the
virtual localization formula of \cite{grabpand1}.   
In fact, virtual localization determines all gravitational 
descendents of $\P^n$. However, these four methods are computationally 
much more complex than the recursions obtained from (\ref{ssmdifff}).

Our recursions for $\P^2$ are most closely related to the Virasoro
conjecture. In \cite{ex}, the authors use weak topological recursion
relations in genus 2 together with the Virasoro conjecture for $\P^2$
to obtain recursions involving a fixed number of descendent integrals 
in each degree (and Gromov-Witten invariants of lower genus).
The Virasoro conjecture generates enough relations to solve
for these series. The numbers $\N_{d}$ below agree  with the values
predicted in \cite{ex}. If the method of \cite{ex} is coupled with 
Getzler's stronger topological recursion relations (\cite{getzler2}),
then the Virasoro conjecture exactly yields recursions involving only
the 3 basic series (\ref{threeser}). It would be quite interesting to
link our recursions to those predicted by the conjecture.

$$\begin{tabular}{R|L|L|L}
d &\N_d &\HH_d &\PP_d \\ \hline
1 &0                       &0                        &0   \\
2 &0                       &0                        &0    \\
3 &0                       &-1/4                     &-1/12 \\
4 &27                      &42                       &25/4   \\
5 &36855                   &130431                   &21119   \\
6 &58444767                &239431851                &33238513      \\
7 &122824720116            &530315850624             &63738316894    \\
8 &346860150644700         &1532247146604636         &161943939423280  \\
9 &1301798459308709880     &5811753079971551880      &547601957576517600 \\
10&6383405726993645784000  &28632855467501316224640  &2432759415312389538720
\end{tabular}$$

\end{document}